\catcode`[=\active  \catcode`]=\active 


\magnification=\magstep1
\baselineskip =5.5mm 
\lineskiplimit =1.0mm
\lineskip =1.0mm

\long\def\comment#1{}

\let\properlbrack=\lbrack
\let\properrbrack=\rbrack
\def\ordcomma{,}
\def\ordcolon{:}
\def\ordsemicolon{;}
\def\ordleftparen{(}
\def\ordrightparen{)}
\def\ordleftbrack{\properlbrack}
\def\ordrightbrack{\properrbrack}
\def\rmcomma{\ifmmode ,\else \/{\rm ,}\fi}
\def\rmcolon{\ifmmode :\else \/{\rm :}\fi}
\def\rmsemicolon{\ifmmode ;\else \/{\rm ;}\fi}
\def\rmleftparen{\ifmmode (\else \/{\rm (}\fi}
\def\rmrightparen{\ifmmode )\else \/{\rm )}\fi}
\def\rmleftbrack{\ifmmode \properlbrack\else \/{\rm \properlbrack}\fi}
\def\rmrightbrack{\ifmmode \properrbrack\else \/{\rm \properrbrack}\fi}
\catcode`,=\active 
\catcode`:=\active 
\catcode`;=\active 
\catcode`(=\active 
\catcode`)=\active 
\let,=\ordcomma
\let:=\ordcolon
\let;=\ordsemicolon
\let(=\ordleftparen
\let)=\ordrightparen
\let[=\ordleftbrack
\let]=\ordrightbrack
\let\lbrack=\ordleftbrack
\let\rbrack=\ordrightbrack
\def\proclaimsl#1{{\sl
\let,=\rmcomma
\let:=\rmcolon
\let;=\rmsemicolon
\let(=\rmleftparen
\let)=\rmrightparen
\let[=\rmleftbrack
\let]=\rmrightbrack
\let\lbrack=\rmleftbrack
\let\rbrack=\rmrightbrack
#1}}

\def\writemonth#1{\ifcase#1
\or January\or February\or March\or April\or May\or June\or July%
\or August\or September\or October\or November\or December\fi}

\newcount\mins
\newcount\minmodhour
\newcount\hour
\newcount\hourinmin
\newcount\ampm
\newcount\ampminhour
\newcount\hourmodampm
\def\writetime#1{%
\mins=#1%
\hour=\mins \divide\hour by 60
\hourinmin=\hour \multiply\hourinmin by -60
\minmodhour=\mins \advance\minmodhour by \hourinmin
\ampm=\hour \divide\ampm by 12
\ampminhour=\ampm \multiply\ampminhour by -12
\hourmodampm=\hour \advance\hourmodampm by \ampminhour
\ifnum\hourmodampm=0 12\else \number\hourmodampm\fi
:\ifnum\minmodhour<10 0\number\minmodhour\else \number\minmodhour\fi
\ifodd\ampm p.m.\else a.m.\fi
}

\font\sevenrm=cmr7
\font\fiverm=cmr5
\font\sevenex=cmex10 at 7pt
\font\sevenbf=cmbx7
\font\fivebf=cmbx5
\font\sevenit=cmti7
\font\seventeenrm=cmr10 scaled \magstep3
\font\twelverm=cmr10 scaled \magstep2
\font\seventeeni=cmmi10 scaled \magstep3
\font\twelvei=cmmi10 scaled \magstep2
\font\seventeensy=cmsy10 scaled \magstep3
\font\twelvesy=cmsy10 at 12pt
\font\seventeenex=cmex10 scaled \magstep3
\font\seventeenbf=cmbx10 scaled \magstep3
\def\sevensize{\sevenrm \baselineskip=4.5 mm%
\textfont0=\sevenrm \scriptfont0=\fiverm \scriptscriptfont0=\fiverm%
\def\rm{\fam0 \sevenrm}%
\textfont1=\seveni \scriptfont1=\fivei \scriptscriptfont1=\fivei%
\def\mit{\fam1 } \def\oldstyle{\fam1 \seveni}%
\textfont2=\sevensy \scriptfont2=\fivesy \scriptscriptfont2=\fivesy%
\def\cal{\fam2 }%
\textfont3=\sevenex \scriptfont3=\sevenex \scriptscriptfont3=\sevenex%
\def\bf{\fam\bffam\sevenbf} \textfont\bffam\sevenbf
\scriptfont\bffam=\fivebf \scriptscriptfont\bffam=\fivebf
\def\it{\fam\itfam\sevenit} \textfont\itfam\sevenit
}
\def\seventeensize{\seventeenrm \baselineskip=5.5mm%
\textfont0=\seventeenrm \scriptfont0=\twelverm \scriptscriptfont0=\sevenrm%
\def\rm{\fam0 \seventeenrm}%
\textfont1=\seventeeni \scriptfont1=\twelvei \scriptscriptfont1=\seveni%
\def\mit{\fam1 } \def\oldstyle{\fam1 \seventeeni}%
\textfont2=\seventeensy \scriptfont2=\twelvesy \scriptscriptfont2=\sevensy%
\def\cal{\fam2 }%
\textfont3=\seventeenex \scriptfont3=\seventeenex%
\scriptscriptfont3=\seventeenex%
\def\bf{\fam\bffam\seventeenbf} \textfont\bffam\seventeenbf
}

\def\setheadline #1\\ #2 \par{\headline={\ifnum\pageno=1 
\hfil
\else \sevensize \noindent
\ifodd\pageno \hfil #2\hfil \else
\hfil #1\hfil \fi\fi}}

\footline={\ifnum\pageno=1 \fiverm \hfil
This paper is not in final form. Typeset using \TeX\ on
\writemonth\month\ \number\day, \number\year\ at \writetime{\time}\hfil 
\else \rm \hfil \folio \hfil \fi}

\def\finalversion{\footline={\ifnum\pageno=1 \hfil 
\else \rm \hfil \folio \hfil \fi}}

\def\moreproclaim{\par}

\def\Head #1:{\medskip\noindent{\bf #1:}}
\def\Proof:{\medskip\noindent{\bf Proof:}}
\def\Proofof #1:{\medskip\noindent{\bf Proof of #1:}}
\def\endproof{\nobreak\hfill$\sqr$\bigskip\goodbreak}
\def\itemi{\item{i)}}
\def\itemii{\item{ii)}}
\def\itemiii{\item{iii)}}
\def\itemiv{\item{iv)}}
\def\itemv{\item{v)}}
\def\itemvi{\item{vi)}}

\def\Abstract\par#1\par{\centerline{\vtop{
\sevensize
\abovedisplayskip=6pt plus 3pt minus 3pt
\belowdisplayskip=6pt plus 3pt minus 3pt
\moreabstract\parindent=0 true in%
A{\fiverm BSTRACT}: \ \ #1}}
\abovedisplayskip=12pt plus 3pt minus 9pt
\belowdisplayskip=12pt plus 3pt minus 9pt
\vskip 0.4 true in}
\def\moreabstract{%
\par \hsize = 5 true in \hangindent=0 true in \parindent=0.5 true in}

\outer\def\firstbeginsection#1\par{\bigskip\vskip\parskip
\message{#1}\leftline{\bf#1}\nobreak\smallskip\noindent}

\outer\def\proclaim#1. #2\par{
\medbreak
\noindent{\bf #1.\enspace}\proclaimsl{#2}\par
\ifdim\lastskip<\medskipamount \removelastskip
\penalty55\medskip\fi}

\def\sqr{\vcenter {\hrule height.3mm
\hbox {\vrule width.3mm height 2mm \kern2mm
\vrule width.3mm } \hrule height.3mm }}

\def\references#1{{
\frenchspacing
\let,=\rmcomma
\let:=\rmcolon
\let;=\rmsemicolon
\let(=\rmleftparen
\let)=\rmrightparen
\let[=\rmleftbrack
\let]=\rmrightbrack
\let\lbrack=\rmleftbrack
\let\rbrack=\rmrightbrack
\halign{\bf##\hfil & \quad\vtop{\hsize=5.5 true in\parindent=0pt\hangindent=3mm
\strut\rm##\strut\smallskip}\cr#1}}}

\catcode`@=11 
\def\vfootnote#1{\insert\footins\bgroup
\sevensize
\interlinepenalty=\interfootnotelinepenalty
\splittopskip=\ht\strutbox
\splitmaxdepth=\dp\strutbox \floatingpenalty=20000
\leftskip=0pt \rightskip=0pt \spaceskip=0pt \xspaceskip=0pt
\textindent{#1}\footstrut\futurelet\next\fo@t}

\def\footremark{\insert\footins\bgroup
\sevensize\it
\let,=\rmcomma
\let:=\rmcolon
\let;=\rmsemicolon
\let(=\rmleftparen
\let)=\rmrightparen
\let[=\rmleftbrack
\let]=\rmrightbrack
\let\lbrack=\rmleftbrack
\let\rbrack=\rmrightbrack
\interlinepenalty=\interfootnotelinepenalty
\splittopskip=\ht\strutbox
\splitmaxdepth=\dp\strutbox \floatingpenalty=20000
\leftskip=0pt \rightskip=0pt \spaceskip=0pt \xspaceskip=0pt
\noindent\footstrut\futurelet\next\fo@t}
\catcode`@=12

\def\Bbb{\bf}
\def\E{{\Bbb E}}
\def\R{{\Bbb R}}
\def\Z{{\Bbb Z}}
\def\N{{\Bbb N}}
\def\C{{\Bbb C}}
\font\eightrm=cmr10 at 8pt
\def\R{\hbox{\rm I\kern-2pt R}}
\def\Z{\hbox{\rm Z\kern-3pt Z}}
\def\N{\hbox{\rm I\kern-2pt I\kern-3.1pt N}}
\def\C{\hbox{\rm \kern0.7pt\raise0.8pt\hbox{\eightrm I}\kern-4.2pt C}} 
\def\E{\hbox{\rm I\kern-2pt E}}

\def\Id{\hbox{\rm Id}}

\def\invc{{c^{-1}}}

\def\list#1,#2{#1_1$, $#1_2,\ldots,$\ $#1_{#2}}

\def\span#1{\overline{\hbox{\rm span}}\{#1\}}

\def\lnorm{\left\|}
\def\rnorm{\right\|}
\def\normo#1{\lnorm #1 \rnorm}
\def\widedot{\,\cdot\,}

\def\trinormo#1{\left|\left|\left| #1 \right|\right|\right|}
\def\trinormdot{\trinormo{\widedot}}

\def\lmod{\left|}
\def\rmod{\right|}
\def\modo#1{\lmod #1 \rmod}

\def\dom#1{{\vert_{#1}}}


\def\couple#1{(#1_0,#1_1)}
\def\sig#1{#1_0 + #1_1}
\def\del#1{#1_0 \cap #1_1}
\def\Xq{X_0\oplus_Q X_1}


\setheadline SPACES OBTAINED BY INTERPOLATION\\
             GARLING -- MONTGOMERY-SMITH

\finalversion

{
\seventeensize
\centerline{\bf Complemented Subspaces of Spaces}
\centerline{\bf Obtained by Interpolation}
}
\bigskip\bigskip
\centerline{\bf D.J.H.~Garling}
{
\sevensize\baselineskip=4.0mm
\centerline{\it St.~John's College, Cambridge CB2 1TP, England.}
}
\medskip
\centerline{\bf S.J.~Montgomery-Smith%
\footnote{}%
{The second named author was supported in part by N.S.F.\ Grant DMS 9001796.}%
}
{
\sevensize\baselineskip=4.0mm
\centerline{\it Department of Mathematics, University of Missouri,}
\centerline{\it Columbia, MO 65211, U.S.A.}
}
\bigskip

\Abstract

If $Z$\ is a quotient of a subspace of a separable Banach space $X$, and
$V$\ is any separable Banach space, then there is a Banach couple
$\couple A$\ such that $A_0$\ and $A_1$\ are isometric to $X\oplus V$,
and any intermediate space obtained using the real or complex interpolation
method contains a complemented subspace isomorphic to $Z$. Thus many properties
of Banach spaces, including having non-trivial cotype, having the Radon--Nikodym
property, and having the analytic unconditional martingale difference
sequence property, do not pass to intermediate spaces.

\footremark{A.M.S.\ (1980) subject classification: 46B99.}

There are many Banach space properties that pass to spaces
obtained by the complex method of interpolation. For example, it is
known that if a couple $\couple A$\ is such that $A_0$\ and $A_1$\ both have
the UMD (unconditional martingale difference sequence) property, and if
$A_\theta$\ is the space obtained using the complex
interpolation method with parameter $\theta$, then $A_\theta$\
has the UMD property whenever $0<\theta<1$. Another example is type of Banach
spaces: if $A_0$\ has type $p_0$\ and $A_1$\ has type $p_1$, then $A_\theta$\
has type $p_\theta$, where $1/p_\theta =
(1-\theta)/p_0 + \theta/p_1$.

Similar results are true for the real method of interpolation. If we denote by
$A_{\theta,p}$\ the space obtained using the real interpolation method from a
couple $\couple A$\ with parameters $\theta$\ and $p$, then $A_{\theta,p}$\ has
the UMD property whenever $A_0$\ and $A_1$\ have the UMD property, $0<\theta<1$,
and $1<p<\infty$. Similarly, if $A_0$\ has type $p_0$\ and $A_1$\ has type
$p_1$, then $A_{\theta,p}$\ has type $p_\theta$, where $1/p_\theta =
(1-\theta)/p_0 + \theta/p_1$\ and $p=p_\theta$\ (see [{\bf 5}] 2.g.22).

However, there are other properties for which it has been hitherto unknown
whether they pass to the intermediate spaces. Examples include the
Radon--Nikodym property, the AUMD (analytic unconditional martingale difference
sequence) property, and having non-trivial cotype. 

This paper deals with these properties, showing that they do not pass to the
intermediate spaces. Indeed, we show the surprising fact that there is a couple
$\couple A$\ such that $A_0$\ and $A_1$\ are both isometric to $l_1$, but all
the real or complex intermediate spaces contain a complemented subspace
isomorphic to $c_0$. This improves a result due to Pisier, who gave an example
of a couple $\couple A$\ for which $A_0$\ is isometric to $L_1$, $A_1$\ is
isometric to a dense subspace of $c_0$, and $c_0$\ is finitely represented in
every intermediate space $A_\theta$\ obtained by the complex interpolation
method (see [{\bf 3}]).

\beginsection Notation

Here we outline the notation we will use about interpolation couples. The reader
is referred to [{\bf 1}] or [{\bf 2}] for details.

A {\it Banach couple\/} is a pair of Banach spaces $\couple A$\ such that
$A_0$\ and $A_1$\ both embed into a common topological vector space, $\Omega$,
which we shall call the {\it ambient space}. Given such a couple, we define
Banach spaces $\sig A$\ (with norm $\normo x = \inf \{\,\normo {x_0}_{A_0} +
\normo {x_1}_{A_1} : x_0 \in A_0$, $x_1\in A_1$, $x_0+x_1 = x \}$) and $\del A$\
(with norm $\normo x = \max\{\normo x_{A_0},\normo x_{A_1} \}$). A map between
two couples $T:\couple A \to \couple B$\ is a linear map $T:\sig A \to \sig B$\
such that $T(A_0) \subseteq B_0$\ and $T(A_1) \subseteq B_1$, and such that
$\normo T_{A_0 \to B_0}$, $\normo T_{A_1\to B_1} < \infty$.

An {\it interpolation method}, $I$, is a functor that takes a Banach couple
$\couple A$\ to a single Banach space $A_I$, such that $\del A \subseteq A_I
\subseteq \sig A$\ with 
$$ \invc \normo x_{\sig A} \le \normo x_{A_I} \le c \, \normo x_{\del A} ,
   \eqno (1) $$
and so that if $T:\couple A \to \couple B$\ is a map between couples, then
$T(A_I) \subseteq B_I$\ with $\normo T_{A_I\to B_I} < \infty$.

An interpolation method is called {\it exponential\/} with {\it exponent
$\theta$\/} if  $0<\theta<1$, and whenever $T:\couple A \to \couple B$\ is a map
between couples, then
$$ \normo T_{A_I \to B_I} \le c \, \normo T_{A_0 \to B_0}^{1-\theta}
   \normo T_{A_1 \to B_1}^\theta . \eqno (2)$$
An interpolation method is called {\it exact exponential\/} if it is
exponential and $c=1$\ in inequalities $(1)$\ and $(2)$.

The most well known interpolation methods are the real interpolation method,
and the complex interpolation method. They are both
exponential, and the complex interpolation method is exact exponential. Another
interpolation method, parameterized by $0<\theta<1$, is the following: if
$x\in \del A$, then let  
$$ \trinormo x = \inf\left\{\,
   \sum_{i=1}^n \normo{x_i}_{A_0}^{1-\theta} \normo{x_i}_{A_1}^\theta :
   \sum_{i=1}^n x_i = x \right\} .$$
Let $A_{\min,\theta}$\ be the completion of $(\del A, \trinormdot)$. This
interpolation method is exact exponential of exponent $\theta$. Furthermore, if
$I$\ is another exponential method of exponent $\theta$, then $\normo
x_{A_I} \le c \, \normo x_{A_{\min,\theta}}$, with $c=1$\ if $I$\ is exact
exponential. (In fact, this interpolation method gives a space equivalent to
the space $A_{\theta,1}$\ obtained by the real interpolation method.)

We denote by $X\oplus V$\ the
Banach space of ordered pairs $(x,v)$\ with norm $\normo{(x,v)} = \normo x
+\normo v$.

All the results given work for Banach spaces over the real or complex scalars.
(Obviously the complex interpolation method will require the complex scalars.)

\beginsection The Main Result

Here we present the main result of this paper.

\proclaim Theorem 1. Suppose that $Z$\ is a Banach space
and that $Q_i$\ is a quotient mapping onto $Z$\ from a closed
linear subspace $Y_i$\ of a separable Banach space $X_i$, for $i=0,1$, such that
the following dimensional constraints hold:
\item{a)} the dimensions of $Y_0$\ and $Y_1$\ are equal;
\item{b)} the codimensions of $Y_0$\ in $X_0$\ and $Y_1$\ in $X_1$\ are equal;
\item{c)} the dimensions of the kernels of $Q_0$\ and $Q_1$\ are equal.
\moreproclaim\noindent
Suppose further that $0<\epsilon<1$, that $0<\theta_0\le \theta_1<1$\ and that
$V_0$\ and $V_1$\ are any infinite dimensional separable Banach spaces. 
\moreproclaim
Then there is a Banach couple $\couple A$, a subspace $W$\ of $\sig A$, and
linear maps $P:W\to Z$\ and $E:Z\to W$\ with the following properties: 
\itemi
$A_0$\ is isometric to $X_0\oplus V_0$, and $A_1$\ is isometric to $X_1\oplus
V_1$; 
\itemii $P\circ E = \Id_Z$;
\moreproclaim\noindent
if $I$\ is an exponential interpolation method of exponent $\theta$, then
\itemiii $A_I \subseteq W$\ and $\normo P_{A_I\to Z} < \infty$\ (with $\normo
P_{A_I\to Z} < 1+\epsilon/3 $\ if $I$\ is exact and
$\theta_0\le\theta\le\theta_1$);
\itemiv $E(Z)\subseteq A_I$\ and $\normo E_{Z\to A_I} < \infty$\ (with $\normo
E_{Z\to A_I} < 1+\epsilon/3 $\ if $I$\ is exact and
$\theta_0\le\theta\le\theta_1$).
\moreproclaim\noindent
Thus $Z$\ is isomorphic to a complemented subspace of $A_I$. (If $I$\ is
exact and $\theta_0\le\theta\le\theta_1$, then $Z$\ is $(1+\epsilon)$\
isomorphic to a $(1+\epsilon)$-complemented subspace of $A_I$.)

This has the following corollaries.

\proclaim Corollary. Given any separable Banach space $X$, there is a Banach
couple $\couple A$\ such that $A_0$\ and $A_1$\ are isometric to $l_1$, and for
every exponential interpolation method $I$, the intermediate space $A_I$\
contains a complemented subspace isomorphic to $X$.

\proclaim Corollary. Given any separable Banach space $X$, there is a
Banach couple $\couple A$\ such that $A_0$\ and $A_1$\ are isometric to
$C([0,1])$, and for every exponential interpolation method $I$, the intermediate
space $A_I$\ contains a complemented subspace isomorphic to $X$.

\proclaim Corollary. The following properties do not pass to
intermediate spaces by any exponential interpolation methods, in the sense
that there is a couple $\couple A$\ such that $A_0$\ and $A_1$\ both have the
property, but none of the intermediate spaces do:
\itemi the Radon--Nikodym property;
\itemii the analytic Radon-Nikodym property;
\itemiii having non-trivial cotype;
\itemiv the AUMD property;
\itemv having a dual with non-trivial cotype;
\itemvi having a dual with the AUMD property.

Note that parts~(v) and~(vi) of the above corollary follow because $C([0,1])^*$\
is finitely represented in $l_1$.

\beginsection A More Elementary Result

Before presenting the proof of Theorem~1, we first prove a more elementary
result that is, in fact, a corollary of Theorem~1. The reason for this is that
the proof of Theorem~1 involves many technicalities, disguising the essential
idea of the proof, which is very simple.

\proclaim Theorem 2. There is a Banach
couple $\couple A$\ such that $A_0$\ and $A_1$\ are isometric to $l_1$, and for
every exponential interpolation method $I$, the intermediate space $A_I$\
contains a complemented subspace isomorphic to $c_0$.

\Proof: Let $e_n$\ be the unit vectors in $l_1$\ and $c_0$, and let $r_n$\ be
an enumeration of the `corners of the unit cube' in $c_0$, that is, all vectors
of the form $(\pm1,\pm1,\ldots,\pm1,0,0,0,\ldots)$. Let $\epsilon_n$\ be the
sequence of numbers defined by
$$ \epsilon_n = (2+\normo{r_n}_1)^{-n} .$$
(Any sequence of numbers tending sufficiently rapidly to zero will do.)
We form an ambient space $\Omega = c_0 \oplus l_1 \oplus l_1$. We define
subspaces
$$ \eqalignno{
   A_0 
   &= \left\{\,
      \left( \sum_n a_n e_n + b_n r_n,
      \sum_n \epsilon_n b_n e_n,
      \sum_n c_n e_n \right)
      : (a_n),\,(b_n),\,(c_n) \in l_1 \,\right\} \cr
   A_1 
   &= \left\{\,
      \left( \sum_n a_n e_n + c_n r_n,
      \sum_n b_n e_n,
      \sum_n \epsilon_n c_n e_n \right)
      : (a_n),\,(b_n),\,(c_n) \in l_1 \,\right\} ,\cr }$$
with norms
$$ \eqalignno{
   \normo{
      \left( \sum_n a_n e_n + b_n r_n,
      \sum_n \epsilon_n b_n e_n,
      \sum_n c_n e_n \right)
   }_{A_0}
   &= \sum_n \left( \modo{a_n} + \modo{b_n} + \modo{c_n} \right) \cr
   \normo{
      \left( \sum_n a_n e_n + c_n r_n,
      \sum_n b_n e_n,
      \sum_n \epsilon_n c_n e_n \right)
   }_{A_1}
   &= \sum_n \left( \modo{a_n} + \modo{b_n} + \modo{c_n} \right) .\cr }$$
The idea is that the unit balls of $A_0$\ and $A_1$ are `slightly perturbed'
versions of the unit balls of $c_0$, where the perturbations for $A_0$\ and
$A_1$\ go in different directions. These perturbations cause the unit balls of
$A_0$\ and $A_1$\ to be convex hulls of linearly independent vectors, and hence
they are affine to the unit ball of $l_1$. The size of these perturbations is
controlled by the quantities $\epsilon_n$. The vectors that are used to perturb
the unit ball of $A_0$\ are also contained in $A_1$\ with no control on their
size. Similarly, the perturbing vectors of $A_1$\ are contained in $A_0$. Thus,
when we form the intermediate space $A_I$, the perturbations in $A_0$\ are
`swamped out' or `washed away' by the vectors in $A_1$, and similarly for the
perturbations in $A_1$, and we are left with a complemented copy of $c_0$. So
much for the idea---now we give the proof.

We define a projection $P:\Omega \to c_0$\ by
$$ P(x,y,z) = x .$$
It is easy to see that $\normo P_{A_i \to c_0} \le 1$\ for $i=0,1$, and
hence $\normo P_{A_I \to c_0} < \infty$\ for any interpolation method $I$.

We define an embedding $E:c_0 \to \Omega$\ by
$$ E(x) = (x,0,0) .$$
We prove that for every
exponential interpolation method that $E(x) \in A_I$\ ($x\in c_0$), with $\normo
E_{c_0\to A_I}<\infty$. To do this, it is sufficient to show that there is a
constant $c$, depending on the exponential interpolation method $I$\ only, such
that for every $n\ge1$\ we have
$$ \normo{(r_n,0,0)}_{A_I} \le c .$$
We first note the following facts.
$$ 
   \normo{(r_n,\epsilon_n e_n,0)}_{A_0} = \normo{(r_n,0,\epsilon_n e_n)}_{A_1} =
   1 $$
$$
   \normo{(0,0,\epsilon_n e_n)}_{A_0} = \normo{(0,\epsilon_n e_n,0)}_{A_1} =
   \epsilon_n $$
$$
   \normo{(r_n,0,0)}_{A_0} = \normo{(r_n,0,0)}_{A_1} = \normo{r_n}_1 .$$
Let us suppose that $I$\ is of exponent $\theta$. We shall give the details in
the case where $I$\ is exact exponential: a similar argument works in the
general case, with less control of the constants. We have the following
inequalities. 
$$ \normo{(r_n,0,0)}_{A_I}
   \le \normo{(r_n,\epsilon_n e_n,\epsilon_n e_n)}_{A_I} 
         + \normo{(0,\epsilon_n e_n,0)}_{A_I} 
         + \normo{(0,0,\epsilon_n e_n)}_{A_I} ,$$
and
$$ \eqalignno{
   \normo{(r_n,\epsilon_n e_n,\epsilon_n e_n)}_{A_I}
   \le & \normo{(r_n,\epsilon_n e_n,\epsilon_n e_n)}_{A_0}^{1-\theta}
         \normo{(r_n,\epsilon_n e_n,\epsilon_n e_n)}_{A_1}^{\theta} \cr
   \le & \biggr(\normo{(r_n,\epsilon_n e_n,0)}_{A_0}
               +\normo{(0,0,\epsilon_n e_n)}_{A_0} \biggl)^{1-\theta} \cr
       & \times \biggr(\normo{(r_n,0,\epsilon_n e_n)}_{A_1}
               +\normo{(0,\epsilon e_n,0)}_{A_1} \biggl)^{\theta} \cr
   \le & (1 + \epsilon_n)^{1-\theta} (1 + \epsilon_n)^\theta \cr
   \le & (1 + \epsilon_n) ,\cr}$$
and
$$ \normo{(0,\epsilon_n e_n,0)}_{A_I}
   \le \normo{(0,\epsilon_n e_n,0)}_{A_0}^{1-\theta}
         \normo{(0,\epsilon_n e_n,0)}_{A_1}^{\theta} .$$
But
$$ \eqalignno{
   \normo{(0,\epsilon_n e_n,0)}_{A_0}
   \le & \normo{(r_n,\epsilon_n e_n,0)}_{A_0} 
         + \normo{(r_n,0,0)}_{A_0} \cr
   \le & 1 + \normo{r_n}_1 ,\cr}$$
and so
$$ \normo{(0,\epsilon_n e_n,0)}_{A_I} \le (1+\normo{r_n}_1)^{1-\theta} \,
   \epsilon_n^{\theta} .$$
Similarly
$$ \normo{(0,0,\epsilon_n e_n)}_{A_I} \le (1+\normo{r_n}_1)^{\theta} \,
   \epsilon_n^{1-\theta} .$$
Therefore
$$ \normo{(r_n,0,0)}_{A_I} \le 1 + \epsilon_n +
   (1+\normo{r_n}_1)^{1-\theta}\epsilon_n^{\theta}
   + (1+\normo{r_n}_1)^{\theta}\epsilon_n^{1-\theta} .$$
By our choice of $\epsilon_n$, this is bounded by some constant $c$\ that
depends only upon $\theta$.

Thus we have bounded maps
$$ P:A_I \to c_0 \qquad E:c_0 \to A_I $$
such that $P \circ I = \Id_{c_0}$, and hence $c_0$ is isomorphic to a
complemented subspace of $A_I$.

Finally, it is very
clear that $A_0$\ and $A_1$\ are both isometric to $l_1$.
\endproof

\beginsection The Proof of Theorem~1

In the sequel, we will make much use of
biorthogonal systems. A {\it biorthogonal system\/} of a Banach space $X$\ is a
sequence $(x_n;\xi_n) \in X\times X^*$\ such that $\xi_n(x_m) = 0$\ if $n\ne m$\
and $1$\ if $n=m$. A biorthogonal system is called {\it
fundamental\/} if $\span{x_n} = X$, and it is called {\it total\/} if $x=0$\
whenever $\xi_n(x) = 0$\ for all $n$. (Of course, if $X$\ is finite
dimensional, then the sequences will be finite.) A result of
Markushevich shows that every separable Banach space has a total, fundamental
biorthogonal system (see [{\bf 6}] or [{\bf 4}] 1.f.3).

We shall need the following proposition and its corollary.

\proclaim Lemma 3. Let $Z$\ be a quotient of a separable
Banach space $Y$\ by the quotient map $Q$\ with 
kernel $K$. Suppose that $(k_p;\kappa'_p)$\ is a total fundamental
biorthogonal system for $K$, and that $(z_n;\zeta_n)$\ is a total fundamental
biorthogonal system for $Z$. Then there are sequences $(y_n)$\ in $Y$, and
$(\psi_n)$\ and $(\kappa_p)$\ in $Y^*$, for which the following hold:
\itemi $(y_n,k_p;\psi_n,\kappa_p)$\ is a total fundamental biorthogonal system
for $Y$;
\itemii $Q(y_n) = z_n$\ and $\zeta_n\circ Q=\psi_n$\ for all $n$;
\itemiii $y\in K$\ if and only if $\psi_n(y) = 0$\ for all $n$.

\Proof: Choose $\tilde y_n \in Y$\ such that $Q(\tilde y_n)
= z_n$, and let $\psi_n = \zeta_n\circ Q$. Using the Hahn--Banach Theorem,
extend $\kappa'_p$\ to $\kappa_p$\ on $Y$\ so that $\kappa_p(\tilde y_n) = 0$\
for $n<p$. Set 
$$ y_n = \tilde y_n - \sum_{m=1}^n \kappa_m(\tilde y_n) k_m . $$
It is now straightforward to verify the conclusions of the lemma.
\endproof

\proclaim Lemma 4. Let $Y$\ be a closed subspace of a separable Banach space
$X$, and let $(y_n;\psi'_n)$\ be a total fundamental biorthogonal system for
$Y$. Then there are sequences $(x_m)$\ in $X$, and $(\xi_m)$\ and $(\psi_n)$\
in $X^*$, such that $(x_m,y_n;\xi_m,\psi_n)$\ is a total fundamental biorthogonal
system for $X$, and such that an element $x$\ in $X$\ belongs to $Y$\ if and
only if $\xi_m(x) = 0$\ for all $m$.

\Proof: By Markushevich's result, $Y$\ has a total fundamental system 
$(y_n;\psi'_n)$, and $X/Y$\ has a total fundamental system $(z_m;\zeta_m)$.
Then the result follows by applying Lemma~3 to the quotient mapping $Q:X\to
X/Y$. 
\endproof

\Proofof Theorem 1: We begin by defining some biorthogonal systems. First, let
$(z_n,\zeta_n)$\ be a total fundamental biorthogonal system for $Z$, with
$\normo{\zeta_n} = 1$\ for all $n$. Using Lemmas~3 and~4, for $i=0,1$, we can
find total fundamental biorthogonal systems
$$ (x_m^i,y_n^i,k_p^i;\xi_m^i,\psi_n^i,\kappa_p^i) $$
for $X_i$\ with the following properties:
\itemi $Q_i(y_n^i) = z_n$\ and $\zeta_n\circ Q_i = \psi_n^i$\ for all $n$;
\itemii $Y_i = \{\,x\in X_i : \xi_m^i(x) = 0 \hbox{ for all }m \}$;
\itemiii $(y_n^i,k_p^i;\psi_n^i\dom{Y_i},\kappa_p^i\dom{Y_i})$\ is a total
fundamental biorthogonal system for $Y_i$;
\itemiv $(k_p^i;\kappa_p^i\dom{K_i})$\ is a total fundamental biorthogonal
system for $K_i$, the kernel of $Q_i$;
\itemv $\normo{\xi_m^i} = \normo{\psi_n^i} = \normo{\kappa_p^i} = 1$\ for all
$m$, $n$\ and $p$.

\noindent
We also choose total fundamental biorthogonal systems $(v_q^i;\phi_q^i)$\ for
$V_i$, with $\normo{\phi_q^i}=1$\ for all $q$. If any of these systems is
finite, then we extend it to an infinite system by including zero terms.

Next we define some sequences of rapidly decreasing numbers. Let $\nu_n$\ and
$\mu_n$\ be sequences of positive numbers for which
$$ \eqalignno{
   \sum L_n \nu_n &\le \epsilon/12 ,\cr
   \max\{\mu_n^\theta L_n,\mu_n^{1-\theta} L_n\} &\le \nu_n 
   \quad\qquad \hbox{for $\theta_0\le\theta\le\theta_1$} ,\cr
\noalign{\hbox{and}}
   \sum \mu_n^\theta L_n &< \infty
   \quad\qquad \hbox{for $0<\theta<1$} ,\cr }$$
where
$$ L_n = \max\{\,
   \normo{x_n^i}, \normo{y_n^i}, \normo{k_p^i}, \normo{v_{2n-1}^i},
   \normo{v_{2n}^i},1 : i=0,1 \} .$$

Now we define the Banach couple. We take as the ambient space $\Omega =
l_\infty\oplus l_\infty\oplus l_\infty\oplus l_\infty\oplus l_\infty $. We
define linear maps $M_i:X_i\oplus V_i\to \Omega$\ in the following way:
$$ \eqalignno{
   M_0(x,v) &= \Bigl(
   \bigl(\psi_n^0(x)\bigr) , 
   \bigl(\xi_m^0(x)\bigr) , \bigl(\mu_p\kappa_p^0(x)\bigr) , 
   \bigl(\mu_q\phi_{2q-1}^0(v)\bigr) , \bigl(\phi_{2q}^0(v)\bigr) 
   \Bigr) ,
   \cr
   M_1(x,v) &= \Bigl(
   \bigl(\psi_n^1(x)\bigr) , 
   \bigl(\mu_q\phi_{2q-1}^1(v)\bigr) , \bigl(\phi_{2q}^1(v)\bigr) ,
   \bigl(\xi_m^1(x)\bigr) , \bigl(\mu_p\kappa_p^1(x)\bigr) , 
   \Bigr) .
   \cr }$$
The important feature here is the interchange in the order of the terms. We
note that $M_0$\ and $M_1$\ are one-one maps into $\Omega$. For $i=0,1$, we set
$A_i = M_i(X_i\oplus V_i)$\ with their norms inherited from the domains. We take
$W$\ to be the linear span of all the exponential interpolation spaces $A_I$.

We shall now go through the details for the case when $I$\ is exact and
$\theta_0\le\theta\le\theta_1$. Otherwise the arguments are similar, with less
control of the constants. Let $e_n$\ denote the $n$th unit vector in
$l_\infty$. Then
$$ (0,\mu_n e_n,0,0,0) = \mu_n M_0(x_n^0) = M_1(v_{2n-1}^1) ,$$
so that
$$ \eqalignno{
   \normo{(0,\mu_n e_n,0,0,0)}_{A_0} &\le \mu_n L_n ,\cr
\noalign{\hbox{and}}
   \normo{(0,\mu_n e_n,0,0,0)}_{A_1} &\le L_n .\cr
\noalign{\hbox{Hence}}
   \normo{(0,\mu_n e_n,0,0,0)}_{A_I} &\le \mu_n^{1-\theta} L_n \le \nu_n .\cr}$$
Similarly,
$$ \eqalignno{
   \normo{(0,0,\mu_n e_n,0,0)}_{A_I} &\le \nu_n ,\cr
   \normo{(0,0,0,\mu_n e_n,0)}_{A_I} &\le \nu_n ,\cr
\noalign{\hbox{and}}
   \normo{(0,0,0,0,\mu_n e_n)}_{A_I} &\le \nu_n .\cr}$$

We define linear functionals on $\Omega$\ as follows: if $t=(h,a,b,c,d) \in
\Omega$, then let $\eta_n(t) = h_n$, $\alpha_n(t) = a_n$, $\beta_n(t) = b_n$,
$\gamma_n(t) = c_n$\ and $\delta_n(t) = d_n$. Then we
have that $$ \eqalignno{
   \normo{\alpha_n}_{A_0^*} &= 1 \cr
\noalign{\hbox{and}}
   \normo{\alpha_n}_{A_1^*} &= \mu_n , \cr
\noalign{\hbox{and so}}
   \normo{\alpha_n}_{A_I^*} &\le \mu_n^\theta \le \nu_n . \cr }$$
Similarly, $\normo{\beta_n}_{A_I^*} \le \nu_n$,
$\normo{\gamma_n}_{A_I^*} \le \nu_n$\ and $\normo{\delta_n}_{A_I^*} \le \nu_n$.

Now we define $E:Z\to \Omega$\ by $E(z) =
\Bigl(\bigl(\zeta_n(z)\bigr),0,0,0,0\Bigr)$. We assert that $E(z) \in A_I$, with
$\normo{E(z)}_{A_I} \le (1+\epsilon/3) \normo z$. To show this, we can
suppose that $\normo z<1$. For $i=0,1$, choose $y_i\in Y_i$\ such that $Q_i(y_i)
= z$\ and $\normo{y_i}\le 1$. Let
$$ v_0 = \sum_p \mu_p \kappa_p^1(y_1) v_{2p}^0 \quad\hbox{and}\quad
   v_1 = \sum_p \mu_p \kappa_p^0(y_0) v_{2p}^1 . $$
It follows from the choice of $\mu_p$\ that these sums converge in $V_0$\ and
$V_1$, and that $\normo{v_0} \le \epsilon/12$\ and $\normo{v_1} \le \epsilon/12$.
Now $$ \eqalignno{
   M_0(y_0,v_0) &= M_1(y_1,v_1) \cr
   &= \Bigl(\bigl(\zeta_n(z)\bigr) , 0 , \bigl(\mu_p\kappa_p^0(y_0)\bigr) ,
      0 , \bigl(\mu_p\kappa_p^1(y_1)\bigr) \Bigr) \cr
   &= w , \cr}$$
say, so that $w\in A_I$\ with $\normo w_{A_I} \le 1+\epsilon/12$. But
$$ \normo{\Bigl( 0 , 0 , \bigl(\mu_p\kappa_p^0(y_0)\bigr) , 0 , 0 \Bigr)}_{A_I}
   \le \sum \mu_p L_p \le \epsilon/12 ,$$
and similarly,
$$ \normo{\Bigl( 0 , 0 , 0 , 0 , \bigl(\mu_p\kappa_p^0(y_0)\bigr) \Bigr)}_{A_I}
   \le \epsilon/12 .$$
Hence $E(z) \in A_I$\ with $\normo{E(z)} \le 1+\epsilon/3$, as desired.

Now we turn to the problem of defining $P$. To do this, we will first show how
we may consider $Z$\ as a subspace of what we shall call
$\Xq$. We let 
$$ D = \{\, (y_0,y_1) : y_i \in Y_i, Q_0(y_0) + Q_1(y_1) = 0 \} ,$$  
and let $\Xq = X_0\oplus X_1 / D$.
We denote by $Q_D$\ the quotient map from $X_0\oplus X_1$\ to $\Xq$.
We define $J:Z\to\Xq$\ by setting $J\bigl(Q_0(y)\bigr)=Q_D(y,0)$\ for $y\in
Y_0$. It is easy to see that $J$\ is well defined, and
that $J\bigl(Q_1(y)\bigr)=Q_D(0,y)$\ for $y\in Y_1$.

Suppose that $w=Q_D(x_0,x_1) \in \Xq$. Then it is easy to verify
that $\pi_m^0(w) = \xi_m^0(x_0)$\ and $\pi_m^1(w) = \xi_m^1(x_1)$ are
well defined maps. It is also easy to check that $Z$\ is isometric to
the space
$$ J(Z) = \{\,w: \pi_m^0(w) = \pi_m^1(w) =0 \} .$$

Now, we will define two bounded operators $S$\ and $T$\ from $A_I$\ to $\Xq$.
First we will concentrate on $S$, which will, in fact, be defined on $\sig A$.

If $t = M_0(x,v) \in A_0$, we let 
$$ R_0(t) = \sum \mu_m \phi_{2m-1}^0(v) x_m^1 .$$
By the choice of $\mu_m$, the sum converges in $X_1$\ with $\normo{R_0(t)} \le
\normo v$. Thus if we set $S_0(t) = Q_D\bigl(x,R_0(t)\bigr)$, then $S_0$\ is
a norm decreasing map from $A_0$\ to $\Xq$. We define $S_1$\ in a similar way.

Now, if $ t = M_0(x_0,v_0) = M_1(x_1,v_1) \in \del A $,
then
$$ \xi_m^0(x_0) = \mu_m \phi_{2m-1}^1(v_1) = \xi_m^0\bigl(R_1(t)\bigr) $$
and
$$ \xi_m^1(x_1) = \mu_m \phi_{2m-1}^0(v_0) = \xi_m^1\bigl(R_0(t)\bigr) ,$$
so that $x_0-R_1(t) \in Y_0$\ and $x_1-R_0(t) \in Y_1$. Further,
$$ \zeta_n\bigl(Q_0(x_0-R_1(t))\bigr) = \psi_n^0(x_0-R_1(t)) = \psi_n^0(x_0),$$
and similarly
$$ \zeta_n\bigl(Q_1(x_1-R_0(t))\bigr) = \psi_n^1(x_1) .$$
Then, since $\psi_n^0(x_0) = \psi_n^1(x_1) = \eta_n(t)$, and since
$(\zeta_n)$\ is total, it follows that $Q_0\bigl(x_0-R_1(t)\bigr) =
Q_1\bigl(x_1-R_0(t)\bigr) $. Thus we see that
$$ S_0(t) = Q_D\bigl(x_0-R_1(t)\bigr) = Q_D\bigl(x_1-R_0(t)\bigr) = S_1(t) ,$$
and so we may define $S:\sig A \to \Xq$\ by setting $S\dom{X_0} = S_0$\
and $S\dom{X_1} = S_1$. Clearly, $S$\ maps $A_I$\
into $\Xq$\ with $\normo S_{A_I\to\Xq} \le 1$.

Now we establish the effect of $\pi_m^0$\ and $\pi_m^1$\ on $S$. If
$t=(h,a,b,c,d) = M_0(x,v) \in A_0$, then
$$ \pi^0_m\bigl(S_0(t)\bigr)
   = \pi_m^0 \bigl(Q_D(x-R_0(t))\bigr)
   = \xi_m^0(x)
   = a_m ,$$
and
$$ \pi^1_m\bigl(S_0(t)\bigr)
   = \pi_m^0 \bigl(Q_D(x-R_0(t))\bigr)
   = \xi_m^0\bigl(R_0(t)\bigr)
   = \mu_m \phi_{2m-1} (v) = c_m .$$
Similarly, if $t=(h,a,b,c,d) = M_1(x,v) \in A_1$, then
$\pi^0_m\bigl(S_1(t)\bigr) = a_m$\ and $\pi^1_m\bigl(S_1(t)\bigr) =
c_m$. Thus, if $t=(h,a,b,c,d) \in A_I$, then
$$ \pi^0_m\bigl(S(t)\bigr) = a_m
   \quad\hbox{and}\quad
   \pi^1_m\bigl(S(t)\bigr) = c_m .$$

Next, we define $T:A_I \to \Xq$\ by
$$ T(t) = Q_D\left(\sum \alpha_m(t) x_m^0 , \sum \gamma_m(t) x_m^1 \right) .$$
Note that the sums converge, and that
$$ \normo T \le \sum \nu_m(\normo{x_m^0} + \normo{x_m^1}) \le \epsilon/6 .$$
We see that, if $t=(h,a,b,c,d)$, then $\pi^0_m\bigl(T(t)\bigr) =
a_m$\ and $\pi^1_m\bigl(T(t)\bigr) = c_m$. 

Hence $(S-T)(t) \in J(Z)$. So now we
set $P=J^{-1}\circ(S-T)$. We have that $\normo{P}_{A_I\to Z} \le 1+\epsilon/3$,
and
$$ \eqalignno{
   P\circ E(z_n) 
   &= P(u_n,0,0,0,0) \cr
   &= J^{-1}\circ(S-T)(u_n,0,0,0,0) \cr
   &= J^{-1}\circ S(u_n,0,0,0,0) \cr
   &= J^{-1}\circ S_0 \circ M_0(y_n^0,0) \cr
   &= J^{-1}\circ Q_D (y_n^0,0) = z_n ,\cr }$$
so that $P\circ E = \Id_Z$.
\endproof

\beginsection Further Results and Conjectures

An obvious extension to Theorem~1 would be the following.

\proclaim Conjecture. Suppose that $\couple Z$\ is a Banach couple where $Z_i$\
is a quotient of a subspace of $X_i$, for $i=0,1$, suppose that $X_0$, $X_1$,
$V_0$\ and $V_1$\ are separable Banach spaces, and suppose that suitable
dimension constraints are satisfied. Then there is a Banach couple $\couple A$\
such that 
\itemi $A_0$\ is isometric to $X_0\oplus V_0$, and $A_1$\ is isometric
to $X_1\oplus V_1$; 
\itemii for every exponential interpolation method $I$, the
space $A_I$\ contains a complemented subspace isomorphic to $Z_I$.
\moreproclaim

The second named author has a tentative proof of a local version of this result
for the real interpolation method, which he hopes to publish later.

However, we can show the following.

\proclaim Theorem 5. Suppose $\couple Y$\ is a Banach couple, where $Y_1
\subset Y_0$\ with $\normo y_{Y_0} \le c \normo y_{Y_1}$. Suppose further that
$Y_0$\ is a quotient space of a separable Banach space $X_0$, and $Y_1$\ is a
complemented subspace of a separable Banach space $X_1$, such that the
codimensions of $Y_0$\ in $X_0$\ and $Y_1$\ in $X_1$\ are equal. Then there is
a Banach couple $\couple A$\ such that 
\itemi $A_0$\ is isometric to $X_0$, and $A_1$\ is isometric to $X_1$;
\itemii given any exponential interpolation method of exponent $\theta$, there
are maps $E:Y_{\min,\theta} \to A_I$\ and  $P:A_I \to Y_I$, both of bounded
norm, such that $P\circ E = \Id_{Y_{\min,\theta}}$
\moreproclaim

\proclaim Corollary. There is a Banach couple $\couple A$\ such that $A_0$\ is
isometric to $l_1$, and $A_1$\ is isometric to $l_p$, and for any exponential
interpolation method $I$\ of exponent $\theta$, the intermediate space
$A_I$\ contains a subspace $V$ such that $l_{p/\theta,1} \subseteq V \subseteq
l_{p/\theta,\infty}$\ with $\invc \normo x_{p/\theta,\infty} \le \normo x_V \le
c \, \normo x_{p/\theta,1}$.

\proclaim Corollary. There is a Banach couple $\couple A$\ such that $A_0$\ has
cotype~$2$, and $A_1$\ has cotype~$p$\ and is $K$-convex, and such that for
all $0<\theta<1$, the real interpolation space $A_{\theta,1}$\ does not have
cotype~$r$\ for any $r < p/\theta$.

This shows that a result of Xu [{\bf 7}] cannot be improved. He showed that, if
$\couple A$\ is a Banach couple such that $A_1$\ is $K$-convex with cotype $p$,
then for all $0<\theta<1$, the real interpolation space $A_{\theta,1}$\ has
cotype $p/\theta$.

\Proofof Theorem 5:  Suppose that $Y_0$\ is a quotient of $X_0$\ by a subspace
$Z_0$, and that $X_1 = Y_1 \oplus Z_1$.  Denote the quotient map from $X_0$\
to $Y_0$\ by $Q$. Let $(z_n^0;\tilde\zeta_n^0)$\ be a total fundamental
biorthogonal sequence for $Z_0$, and let $(z_n^1;\tilde\zeta_n^1)$\ be a total
fundamental biorthogonal sequence for $Z_1$. Extend $\tilde\zeta_n^0$\ to
$\zeta_n^0$\ on $X_0$\ using the Hahn--Banach Theorem, and define
$\zeta_n^1$\ on $X_1$\ by setting $\zeta_n^1(y,z) = \tilde\zeta_n^1(z)$.
Without loss of generality, $\normo{\zeta_n^0} = \normo{\zeta_n^1} = 1$. 

To save chasing constants, let us assume that $I$\ is exact. Let $\epsilon_n$\
be a sequence of sufficiently small positive numbers, so that
$$ \delta = \sum_{n=1}^\infty  
   \max\{\epsilon_n,\epsilon_n^\theta\}
   \max\{\normo{z_n^0},\normo{z_n^1}\}
   < \infty ,$$ 
for all $0<\theta<1$.
We will let $Y_0\oplus
l_\infty $\ be the ambient space. Let $M_i: X_i \to Y_0\oplus l_\infty$\ be the
mappings: 
$$ \eqalignno{
   M_0 (x) &= \Bigl( Q(x), \bigl(\epsilon_n \zeta_n^0(x)\bigr) \Bigr) \cr
   M_1 (x) &= \Bigl( Q(x), \bigl(\zeta_n^1(x)\bigr)  \Bigr) .\cr }$$
Since $(z_n^0;\tilde\zeta_n^0)$\ and $(z_n^1;\tilde\zeta_n^1)$\ are both total
biorthogonal systems, $M_0$\ and $M_1$\ are one to one. We let $A_i =
M_i (X_i)$\ with norms inherited from the domains of
the functions.

Define $E:Y_0 \to Y_0\oplus l_\infty$\ by
$$ E(y) = (y,0) .$$
We desire to show that $E(Y_{\min,\theta}) \subseteq A_I$, with
$\normo{E}_{Y_{\min,\theta}\to A_I} <\infty$. First, note that
$$ \normo{(0,\epsilon_n e_n)}_{A_0} = \normo{z_n^0} ,$$
and
$$ \normo{(0,\epsilon_n e_n)}_{A_1} = \epsilon_n \normo{z_n^1} .$$
Hence, 
$$ \normo{(0,\epsilon_n)}_{A_I} 
   \le \epsilon_n^\theta \max\{\normo{z_n^0},\normo{z_n^1}\} ,$$
and so
$$ \normo{\Bigl(0,\bigl(\epsilon_n \zeta_n^0(x)\bigr)\Bigr)}_{A_I} 
   \le \delta \normo x_{X_0} .$$
Let $\nu > 0$.
For $y \in Y_1$, let $x \in X_0$\ be such that $Q(x) = y$\
with $\normo x_{X_0} \le (1+\nu) \normo y_{Y_0}$. Then
$$ \eqalignno{
   \normo{\Bigl(y,\bigl(\epsilon_n \zeta_n^0(x) \bigr)\Bigr)}_{A_0} 
   &= \normo{x}_{X_0} \le (1+\nu) \normo y_{Y_0} ,\cr
   \normo{(y,0)}_{A_1} &= \normo y_{Y_1} , \cr}$$
and
$$ \eqalignno{
   \normo{\Bigl(0,\bigl(\epsilon_n \zeta_n^0(x) \bigr)\Bigr)}_{A_1} 
   &\le \sum_n \epsilon_n \normo{z_n^1} \normo x_{X_0} \cr
   &\le (1+\nu) \delta \normo y_{Y_0} 
   \le (1+\nu) \delta \normo y_{Y_1} . \cr}$$
So,
$$ \normo{\Bigl(y,\bigl(\epsilon_n \zeta_n^0(x) \bigr)\Bigr)}_{A_1} \le
   (1+\nu)(1+\delta) \normo y_{Y_1} .$$
Therefore,
$$ \normo{\Bigl(y,\bigl(\epsilon_n \zeta_n^0(x) \bigr)\Bigr)}_{A_I} \le
   (1+\nu)(1+\delta) \normo y_{Y_0}^{1-\theta} \normo y_{Y_1}^\theta .$$
Finally,
$$ \eqalignno{
   \normo{(y,0)}_{A_I} 
   &\le \delta \normo x_{X_0} 
   + (1+\nu) (1+\delta) 
   \normo y_{Y_0}^{1-\theta} \normo y_{Y_1}^\theta \cr
   &\le (1+\nu) (1+2\delta) 
  \normo y_{Y_0}^{1-\theta} \normo y_{Y_1}^\theta . \cr}$$
This is sufficient to show that
$\normo{E(y)}_{A_I} \le (1+\nu)(1+2\delta) \normo{y}_{Y_{\min,\theta}}$\
as desired.

Now we let $P:Y_0\oplus l_\infty \to Y_0$\ be the map
$$ P(y,a) = y .$$
Then $\normo P_{A_0 \to Y_0} = \normo P_{A_1 \to Y_1} = 1$, and hence $\normo
P_{A_I \to Y_I} \le 1$. Clearly $P\circ E = \Id_{Y_{\min,\theta}}$.
\endproof

\beginsection Acknowledgements

We would like to express our warm thanks to G.~Pisier for suggesting this line of
research to us.

\beginsection References

\references{
1 & C.~Bennett and R.~Sharpley,\sl\ Interpolation of Operators,\rm\
Academic Press 1988.\cr
2 & J.~Bergh and J. L\"ofstr\"om,\sl\ Interpolation Spaces,\rm\
Springer-Verlag 1976.\cr
3 & S.J.~Dilworth,\rm\ Complex convexity and the geometry of Banach spaces,\sl\
Math.\ Proc.\ Camb.\ Phil.\ Soc.\ {\bf 99} (1986), 495--506.\cr
4 & J.~Lindenstrauss and L.~Tzafriri,\sl\ Classical Banach Spaces
I---Se\-qu\-ence Spa\-ces,\rm\ Springer-Verlag 1977.\cr
5 & J.~Lindenstrauss and L. Tzafriri,\sl\ Classical Banach Spaces
II---Fu\-nc\-t\-ion Spa\-ces,\rm\ Springer-Verlag 1979.\cr
6 & A.I.~Markushevich,\rm\ On a basis in the wide sense for linear spaces,\sl\
Dokl.\ Akad.\ Nauk {\bf 41} (1943), 241--244.\cr
7 & Q.~Xu,\rm\ Cotype of the spaces $(A_0,A_1)_{\theta 1}$,\sl\ Geometric
Aspects of Functional Analysis, Israel Seminar 1985--6, J.~Lindenstrauss
and V.D.~Milman (Eds.),\rm\ Springer Verlag 1987.\cr
}

\bye